\documentclass{amsart}

\usepackage[english]{babel}
\usepackage{blindtext}
\usepackage[utf8]{inputenc}
\usepackage{fancyhdr}
\usepackage{amsmath}
\usepackage{amssymb}
\usepackage{mathtools}
\usepackage{tikz-cd}
\usepackage{tikz-cd}
\usepackage{enumitem}
\usepackage{latexsym,amsfonts,amssymb,amsthm,amsmath}
\usepackage{graphicx}
\usepackage[colorlinks=true, allcolors=blue]{hyperref}

\usepackage{mathtools}

\def\F{\mathcal F}

\def\E{\mathcal E}

\newtheorem{thm}{Theorem}[section]

\newtheorem{lem}[thm]{Lemma}

\theoremstyle{definition}

\theoremstyle{remark}
\newtheorem*{remark}{Remark}

\title[Sections of $K3$ surfaces and Mercat's conjecture]{Sections of $K3$ surfaces with Picard number two and Mercat's conjecture}
\author{Marian Aprodu, Laura Filimon}
\dedicatory{Dedicated to the memory of Lucian B\u adescu}
\thanks{MA was partly supported by the PNRR	grant CF 44/14.11.2022 \emph{Cohomological Hall algebras of smooth surfaces and applications.}}
\address{MA: Simion Stoilow Institute of Mathematics, P.O. Box 1-764, 014700 Bucharest,
	Romania
	\& Faculty of Mathematics and Computer Science, University of Bucharest, Bucharest,
	Romania}
\email{marian.aprodu@imar.ro \& marian.aprodu@fmi.unibuc.ro}
\address{LF: Faculty of Mathematics and Computer Science, University of Bucharest, Bucharest,
	Romania}
\email{laura.filimon@my.fmi.unibuc.ro}
\keywords{Higher-rank Brill-Noether theory, curves on $K3$ surfaces}
\subjclass{Primary 14H60. Secondary 14H51, 14J60, 14J28}

\begin{document}
	\maketitle
	
	\begin{abstract}
		In \cite[Theorem 1.1]{Farkas-Ortega}, the authors present counterexamples to Mercat's conjecture by restricting to a hyperplane section $C$ some suitable rank-two vector bundles on a $K3$ surface whose Picard group is generated by $C$ and another very ample divisor. We prove that the same bundles produce other counterexamples by restriction to hypersurface sections $C_n\in|nC|$ for all $n\ge 2$. In the process, we compute the Clifford indices of the corresponding hypersurface sections $C_n$, noting their non-generic nature for $n\ge 2$ (refer to Theorem \ref{thm:Cliff}). A key ingredient to prove the (semi)stability of the restricted bundles, Theorem \ref{thm:semistability}, is Green's Explicit $H^0$ Lemma (see \cite[Corollary (4.e.4)]{Green}). In what concerns the (semi)stability, although general restriction theorems such as \cite[Theorem 1.2]{Flenner} or \cite[Theorem 1.1]{Feyzbakhsh} are applicable for sufficiently large, explicit values of $n$, our approach works for all $n\ge 2$. It is also worth noting that our proof deviates slightly from the one presented in \cite[Proposition 3.2]{Farkas-Ortega}. 
		Employing the same strategy leads to an enhancement of the main result of \cite{Sengupta}; refer to Theorem \ref{thm:VeryGeneral} for counterexamples to the conjecture on curves in $|nC|$, where $C$ now acts as a generator of the Picard group.
	\end{abstract}
	
	\section{Introduction}
	
	Mercat's conjecture aims to establish a connection between higher-rank Brill--Noether theory and classical Brill--Noether theory concerning curves. Let $C$ be a smooth curve of genus $g\ge 3$, and consider $\E$ a semistable rank vector bundle on $C$ satisfying $h^1(C,\E) \ge h^0(C,\E) \ge 2r$. The \emph{Clifford index of $\E$} is defined as
	\[
	\gamma(\E) := \mu(\E) - \frac{2h^0(C,\E)}{\mathrm{rk}(\E)}+2\ge 0
	\]
	and $\mathrm{Cliff}_r(C)$, the \emph{$r$th Clifford index of $C$} is the minimum of the Clifford indices of bundles of rank $r$ that can contribute, i.e.
	\[
	\mathrm{Cliff}_r(C):=\mathrm{min}\{\gamma(\E) \ :\ \E \in \mathcal{U}_C(r,d), d \le r(g - 1), h^0(C,\E) \ge 2r\}.
	\] 
	
	In this context, Mercat \cite{Mercat} conjectured that for any $r\ge 1$, we have $\mathrm{Cliff}_r(C)=\mathrm{Cliff}(C)$. It is worth noting that the inequality $\mathrm{Cliff}_r(C)\le\mathrm{Cliff}(C)$ is readily obtained taking direct sums of line bundles $A^{\oplus r}$. Originally, the conjecture was formulated as an explicit upper bound in terms of $\mathrm{Cliff}(C)$ for the number of sections of semistable bundles, \cite[p. 786]{Mercat}. Precisely, the conjectured bound is given by:
	\[
	h^0(C,\E)\le \frac{d}{2}-r\left(\frac{\mathrm{Cliff}(C)}{2}-1\right)
	\]
	for all $\E \in \mathcal{U}_C(r,d)$, with $d \le r(g - 1)$ and $h^0(C,\E) \ge 2r$. In rank two, this inequality simplifies to $h^0(C,\E)\le \frac{d}{2}-\mathrm{Cliff}(C)+2$ for all $\E \in \mathcal{U}_C(2,d)$, with $d \le 2g-2$ and $\E$ having at least 4 independent sections. Note that the number of independent sections is always bounded, \cite[Proposition 3, Proposition 4]{Re}, \cite[Theorem 2.1]{Mercat} etc, but the known general bounds are weaker than those predicted by the conjecture.
	
	While the conjecture has been confirmed in various cases, e.g. in rank two, it holds for arbitrary $k$-gonal curves of genus $g > 2(k-1)(k-2)$, for general curves \cite{Bakker-Farkas}, for general $k$-gonal curves of genus $g > 4k-4$, for plane curves \cite{Lange-Newstead},   \cite{Lange-Newstead13} etc, it fails for large values $k$ of the gonality. Specifically, several counterexamples have been provided by curves on $K3$ surfaces, as seen, for instance in \cite{Farkas-Ortega11}, \cite{Lange-Newstead11}, \cite{Lange-Newstead12}, \cite{Farkas-Ortega}, \cite{Sengupta} (see also \cite{Aprodu-Farkas-Ortega}, \cite{Feyzbakhsh} for higher ranks). A current challenge is to discover additional examples of pairs $(g,k)$ where Mercat's conjecture fails in rank two or to determine whether the existing list of counterexamples is exhaustive. In view of \cite[Section 4]{Lange-Newstead13}, the problem needs to be addressed for curves of Clifford dimension one.
	
	In this short Note, we present a new infinite set of counter-examples for the conjecture. Our methodology also revolves around the utilization of curves on $K3$ surfaces, specifically drawing upon the counterexamples identified in \cite{Farkas-Ortega}, i.e. curves on $K3$ surfaces of Picard number two. A distinctive aspect of our investigation is the transition from hyperplane sections to hypersurface sections, aligning with the exploration of $K3$ surfaces with Picard number one, as discussed in \cite{Sengupta}. The two main technical difficulties that we have to overcome are: the computation of Clifford indices, and the semistability of the restricted bundles. These issues are addressed in Theorem \ref{thm:Cliff} and Theorem \ref{thm:semistability}, respectively. For the computation of the Clifford indices, we use the Main Theorem of \cite{Green-Lazarsfeld}, and the verification of semistability relies on Green's explicit $H^0$ Lemma, see \cite[Corollary (4.e.4)]{Green}. Employing the same strategy leads to an enhancement of the main result of \cite{Sengupta}; refer to Theorem \ref{thm:VeryGeneral} for counterexamples to the conjecture on curves in $|nC|$, where $C$ now acts as a generator of the Picard group.

	\section{Basic properties of Lazarsfeld--Mukai bundles}
	\label{sec:preliminaries}

	We follow closely the presentation of \cite{Lazarsfeld}. Let $S$ be a $K3$ surface, $C$ be a smooth connected curve of genus $g$ in $S$, and $A$ be a base-point-free complete $g^r_d$ on $C$. Denote by $M_A$ the kernel of the evaluation map
	\[
	\mathrm{ev}_A:H^0(A)\otimes \mathcal O_C\to A.
	\] 
	
	The map $\mathrm{ev}_A$ induces a surjective morphism $H^0(A)\otimes \mathcal O_S\to A$ of sheaves on $S$
	whose kernel $\F_{C,A}$ is a vector bundle of rank $(r+1)$. Its dual $\E_{C,A}=\F_{C,A}^\vee$ is called a \emph{Lazarsfeld--Mukai bundle}. The defining sequences of $\F_{C,A}$ and $\E_{C,A}$ are
	\begin{equation}
		\label{eqn:F}
		0\to \F_{C,A}\to H^0(A)\otimes \mathcal O_S\to A\to 0.
	\end{equation}	
	and, respectively 
	\begin{equation}
		\label{eqn:E}
		0\to H^0(A)^\vee\otimes \mathcal O_S\to \E_{C,A}\to K_C(-A)\to 0.
	\end{equation}
	
	The bundles $\E_{C,A}$ and $\F_{C,A}$ have the following properties:
	\begin{enumerate}
		\item
		$\mathrm{det}(\E_{C,A})=\mathcal{O}_S(C)$,
		
		\item
		$c_2(\E_{C,A})=d$,
		
		\item
		$h^0(S,\F_{C,A})=h^1(S,\F_{C,A})=0$, 
		
		\item
		$\chi(S,\F_{C,A})=h^2(S,\F_{C,A})=2(r+1)+g-d-1$,
		
		\item
		$h^0(S,E_{C,A})=r+1+h^0(C,K_C(-A))$,
		
		\item
		$\E_{C,A}$ is generated off the base locus of $|K_C(-A)|$ inside $C$.
	\end{enumerate}
	
	Restricting the sequence (\ref{eqn:F}) to the curve $C$, we obtain a short exact sequence:
	\begin{equation}
		\label{eqn:F|C}
		0\to K_C^\vee(A)\to \F_{C,A}|_C\to M_A\to 0
	\end{equation}
	which implies, twisting by $K_C(-A)$ and using the adjunction formula,
	\begin{equation}
		\label{eqn:F|C al doilea}
		0\to \mathcal{O}_C\to \F_{C,A}\otimes K_C(-A)\to M_A\otimes K_C(-A)\to 0.
	\end{equation}
	
	Note that $H^0(M_A\otimes K_C(-A))=\ker(\mu_{0,A})$, where $\mu_{0,A}:H^0(A)\otimes H^0(K_C(-A))\to H^0(K_C)$ is the Petri map.

	\section{Clifford indices of hypersurface sections of a $K3$ surface with Picard number two}
	\label{sec:Cliff}

	Given integers $p\geq 3$ and $a\geq 2p+3$, let $S$ be a $K3$ surface whose Picard group is generated by two very ample smooth divisors, $\text{Pic}(S)=\langle C, D\rangle $, where $C^2=4a$, $D\cdot C=2a+2p+1$, $D^2=4p+2$. The existence of such surfaces is established through the surjectivity of the period map, as noted in \cite{Farkas-Ortega}. We focus on the embedding $S\subset \mathbb{P}^{2a+1}$ defined by the complete linear system $|C|$. It is worth noting that in \cite{Farkas-Ortega}, the authors consider the surface $S$ as being embedded via the other linear system~$|D|$, denoted by $|H|$ in that context.
	
	For the convenience of the reader, we highlight the following simple fact that was implicitly used in \cite{Farkas-Ortega}.

	\begin{lem}
		\label{lem:elliptic}
		Put $E=C-D$. Then $E^2=0$, $h^0(S,\mathcal{O}_S(E)) = 2$ and $h^1(S,\mathcal{O}_S(E))=0$.
	\end{lem}
	
	\proof
	The numerical data makes it evident that $E^2=0$. Notably, as $(-E)\cdot D=-2a+2p+1< 0$ and $D$ is ample, it implies that $-E$ cannot be effective, and it cannot be zero either. By the Riemann-Roch Theorem, we derive $h^0(S,\mathcal{O}_S(E)) \ge 2$.
	
	Suppose the linear system $|E|$ is base-point-free; in this case, according to \cite[Proposition 2.6]{SaintDonat}, it follows that $E$ is a multiple of a smooth elliptic curve. Since $E$ is a generator of the Picard group, it must be a smooth elliptic curve. Consequently, we have $h^0(S,\mathcal{O}_S(E)) = 2$ and $h^1(S,\mathcal{O}_S(E))=0$.

	Now, assume the linear system $|E|$ has base points.  According to \cite[Proposition 2.6]{SaintDonat} the linear system $|E|$ has a fixed component $\Delta$. 
	Write $E=\Delta+E'$ where $E'$ is an effective divisor with $h^0(S,\mathcal{O}_S(E))=h^0(S,\mathcal{O}_S(E'))$. From \cite[Proposition 2.2]{Farkas}, we deduce that $(E')^2\ge 0$ and $\Delta^2\ge 0$. Since $|E'|$ is the moving part of the linear system $|E|$, we must also have $E'\cdot \Delta > 0$. On the other hand $E^2=0$, which is a contradiction.
	\endproof

	We aim to prove that the Clifford index of any curve in the linear system $|nC|$, for $n\ge 2$, is computed by $\mathcal{O}(E)$.

	\begin{thm}
		\label{thm:Cliff}
		For any $n\ge 2$, and any smooth curve $C_n\in|nC|$, we have $\mathrm{Cliff}(C_n)=n(2a-2p-1)-2$.
	\end{thm}
	
	\proof
	We remark that the genus of $C_n$ is $g(C_n)=2an^2+1$, and the bundle $\mathcal{O}_{C_n}(E)$ contributes to the Clifford index of $C_n$, with its Clifford index strictly smaller than the generic Clifford index $(an^2-1)$. Applying the Main Theorem of \cite{Green-Lazarsfeld}, and  \cite[Lemma 2.2]{Martens} the Clifford index of $C_n$ is computed by the restriction of a line bundle $\mathcal{O}_S(F)\in\mathrm{Pic}(S)$ by the formula 
	\[
	\mathrm{Cliff}(C_n)=\mathrm{Cliff}(\mathcal{O}_{C_n}(F))=F\cdot C_n-F^2-2.
	\]
	To simplify calculations, we work with the basis $\{C,E\}$ of $\mathrm{Pic}(S)$ instead of the original $\{C,D\}$, considering $E^2=0$. Note that $C\cdot E=2a-2p-1>0$. Therefore, expressing $F=sC+tE$ with $s,t\in\mathbb{Z}$, we compute:
	\begin{equation}
		\label{eq:f}
		f(s,t):=\mathrm{Cliff}(\mathcal{O}_{C_n}(sC+tE))=(n-2s)(2a-2p-1)t-4as^2+4ans-2.
	\end{equation}
	The condition $f(s,t)\ge 0$ must be satisfied due to the definition of the Clifford index.
	
	Following the proof of \cite[Theorem 3]{Farkas} and the proof of \cite[Proposition 3.3]{Farkas-Ortega}, we observe that $F$ is subject to the following restrictions:
	\begin{enumerate}
		\item[(i)] $F^2\ge 0$,
		\item[(ii)] $F\cdot D> 2$,
		\item[(iii)] $F\cdot C_n\le g(C_n)-1$,
	\end{enumerate}
	
	Taking into account that $g(C_n)=2an^2+1$, these constraints are translated into the following conditions:
	\begin{enumerate}
		\item[(i)] $s(2as+(2a-2p-1)t)\ge 0$,
		\item[(ii)] $4as+(2a-2p-1)(t-s)>2$,
		\item[(iii)] $4as+(2a-2p-1)t\le 2an$.
	\end{enumerate}	
	
	The objective is to prove that the minimum of $f$, when $s$ and $t$ are integers satisfying conditions (i), (ii), and (iii) is attained at $(0,1)$. Since $E\cdot C_n-E^2-2=n(2a-2p-1)-2$ that would conclude the proof of the theorem.

	We note that $s\ge 0$. Indeed, if $s<0$, then (i) implies that $2as+(2a-2p-1)t\le 0$ and hence we obtain from (ii) that $(2p+1)s>0$ which is a contradiction with the assumption $s<0$. Consequently, condition (i) is reformulated as:
	\begin{enumerate}
		\item[(i)] $2as+(2a-2p-1)t\ge 0$.
	\end{enumerate}
	
	Additionally, we have $s\leq n$, due to (i) and (iii) leading to $2an-2as\geq 0$. 
	
	\medskip
	
	If $s=0$, then $f(0,t)=n(2a-2p-1)t-2$ and the minimal positive value is $f(0,1)=n(2a-2p-1)-2$ which we wanted to prove.
	
	\medskip
	
	We analyze next the case $s\ge 1$. The inequalities (i) and (iii) give the following bounds for $t$:
	\[
	t_{min}:=-\frac{2as}{2a-2p-1}\le t\le t_{max}:=\frac{2a(n-2s)}{2a-2p-1}.
	\]
	
	If $n$ is even and $s=\frac{n}{2}$, we observe that $f\left(\frac{n}{2},t\right)=an^2-2>n(2a-2p-1)-2$ for $n\ge 2$.
	
	If $n>2s$, since the coefficient of $t$ in the expression of $f$ is positive and $s\ge 1$, it follows that
	\[
	f(s,t)\ge f(s,t_{min})=2ans-2>n(2a-2p-1)-2.
	\]
	
	If $n<2s$, since the coefficient of $t$ in the expression of $f$ is negative, it holds that
	\[
	f(s,t)\ge f(s,t_{max})=4as^2-4ans+2an^2-2.
	\]
	
	On the interval $\left[\frac{n}{2}, n\right]$ the degree-two function $g(s):=f(s,t_{max})$ is increasing and hence
	\[
	f(s,t)\ge f(s,t_{max})\ge f\left(\frac{n}{2},t_{max}\right)=an^2-2>n(2a-2p-1)-2
	\]
	for $n\ge 2$. This completes the proof.
	\endproof

	\begin{remark}
		For any integer $n\geq 2$, any $K3$ surface $S$, and any very ample line bundle $\mathcal{O}_S(C)$, consider a smooth curve $C_n$ in the linear system $|nC|$. In this context, the Clifford index of $C_n$ is smaller than the generic value $\left[\frac{g(C_n)-1}{2}\right]$. Specifically, the restriction of the bundle $\mathcal{O}_S(C)$ to $C_n$ contributes to the Clifford index, and upon direct computation, its Clifford index is found to be smaller than the generic value. If $\mathcal{O}_S(C)$ generates the Picard group of $S$, then $\mathrm{Cliff}(C_n)$ is computed by the restriction of $\mathcal{O}_S(C)$. In contrast to the very generic case, the explicit situation presented here yields $\mathrm{Cliff}(\mathcal{O}_{C_n}(C))=4(n-1)a-2>n(2a-2p-2)-2=\mathrm{Cliff}(C_n)$.
	\end{remark}

	\section{New counterexamples to Mercat's conjecture}
	\label{sec:counterexamples}

	We adopt the notation from in the previous sections. Consider a $g^1_{p+2}$ denoted as $A$ on $D$, and let $\E=\E_{C,A}$ be the associated Lazarsfeld--Mukai bundle. As affirmed by \cite[Theorem~1.1]{Farkas-Ortega}, it follows that $\text{Cliff}(C)=a$, and additionally, $\gamma(\E|_C)< \text{Cliff}(C)$. Notably, $\E|_C$ is semistable (\cite[Proposition 3.2]{Farkas-Ortega}), consequently providing a counterexample to Mercat's conjecture.
	
	In the subsequent discussion, we establish the following result.
	
	\begin{thm}
		\label{thm:semistability}
		Notation as above. Assume $a\ge 3p+2$. For any $n\ge 2$, the bundle $\E|_{C_n}$ is stable and it is a counter-example to Mercat's conjecture.
	\end{thm}
	
	\proof
	We prove first the semistability of $\E|_{C_n}$. Suppose, for a contradiction, that $\E|_{C_n}$ is not stable and consider
	\[
	0\to \mathcal{O}_{C_n}(B)\to \E|_{C_n}\to \mathcal{O}_{C_n}(D-B)\to 0
	\]
	a destabilizing sequence. In particular,
	\begin{equation}
		\label{eq:deg(B)}
		\textrm{deg}(B) \ge \mu(\E|_{C_n})=\frac{n(2a+2p+1)}{2}.
	\end{equation}
	
	Since $\E$ is globally generated, it follows that $\mathcal{O}_{C_n}(D-B)$, along with any other quotient of $\E$, is also globally generated.
	
	If  $\mathcal{O}_{C_n}(D-B)\ne \mathcal{O}_{C_n}$, then  $h^0(C_n,\mathcal{O}_{C_n}(D-B))\ge 2$. 
	Furthermore, since $h^0(C_n,\mathcal{O}_{C_n}(C_n-D+B))\ge h^0(C_n,\mathcal{O}_{C_n}(C_n-D))\ge h^0(S,\mathcal{O}_S(C_n-D))\ge 2$, it follows that $\mathcal{O}_{C_n}(D-B)$ contributes to the Clifford index of $C_n$. 
	Using the inequality (\refeq{eq:deg(B)}) we evaluate
	\[
	\mathrm{Cliff}(\mathcal{O}_{C_n}(D-B))\le n(2a+2p+1)-\mathrm{deg}(B)-2\le \frac{n(2a+2p+1)}{2}-2
	\]
	and the latter value is smaller than $n(2a-2p-1)-2=\mathrm{Cliff}(C_n)$, by the assumption $a\ge 3p+2$, leading to a contradiction. 
	
	In conclusion, we have $\mathcal{O}_{C_n}(D-B) = \mathcal{O}_{C_n}$. This implies the existence of a short exact sequence
	\[
	0\to \mathcal{O}_{C_n}(D)\to \E|_{C_n}\to \mathcal{O}_{C_n}\to 0
	\]
	and, as a consequence, we have:
	\begin{equation}
		\label{eq:no-sections}
		h^0(C_n, \E|_{C_n})\ge h^0(C_n,\mathcal{O}_{C_n}(D)).
	\end{equation}

	Since $h^0(S,\mathcal{O}_{S}(D-nC))=0$ for all $n\ge 1$, it follows that $h^0(C_n,\mathcal{O}_{C_n}(D)) \ge h^0(S,\mathcal{O}_S(D)) = 2p+1$. Moreover, the two dimensions are equal, as shown below.
	
	\medskip
	
	\emph{Claim 1.} $h^1(S,\mathcal{O}_S(D-nC))=0$ for all $n\ge 1$.
	
	\medskip
	
	We proceed by induction on $n$. For the base case $n=1$, we apply Lemma \ref{lem:elliptic}. For the induction step, for $n\ge 2$, consider the long cohomology sequence of the short exact sequence
	\[
	0\to \mathcal{O}_S((n-1)C-D)\to \mathcal{O}_S(nC-D)\to \mathcal{O}_C(nC-D)\to 0
	\]
	and observe that $h^1(C,K_C^{\otimes n}(-D)))=0$ by degree reasons.
	
	\medskip
	
	\emph{Claim 2.} $h^1(S,\E(-nC))=0$.
	
	\medskip
	
	To establish Claim 2, we begin with the defining exact sequence \eqref{eqn:F} of the dual of the Lazarsfeld-Mukai bundle:	    
	\[
	0\to \E^\vee\to H^0(A)\otimes\mathcal{O}_S\to A\to 0
	\]
	(where $A$ was a $g^1_{p+2}$ on $D$) twist it with $\mathcal{O}_S(nC)$ and take the long cohomology sequence:
	\[
	0\to H^0(\E^\vee(nC))\to H^0(A)\otimes H^0(\mathcal{O}_S(nC))\to H^0(D,A(nC))\to H^1(\E^\vee(nC))\to 0.
	\]
	
	Claim 1 implies that the restriction map $H^0(S,\mathcal{O}_S(nC))\to H^0(D,\mathcal{O}_D(nC))$ is surjective. Hence, Claim 2 would follow from the surjectivity of the multiplication map
	\[
	H^0(D,A)\otimes H^0(D,\mathcal{O}_D(nC))\to H^0(D,A(nC))
	\]
	To this end, we apply \cite[Corollary (4.e.4)]{Green}; the hypothesis 
	$$
	\mathrm{deg}(A)+\mathrm{deg}(\mathcal{O}_D(C_n))\ge 4g(D)+2
	$$ 
	is verified for $n\ge 2$, as the genus of $D$ is $g(D)=2p+2$, and $\mathrm{deg}(A)+\mathrm{deg}(\mathcal{O}_D(C_n))=(p+2)+n(2a+2p+1)\ge 13p+16$. Claim 2 is proved.
	
	\medskip
	
	We consider the short exact sequence
	\[
	0\to \E(-nC)\to \E\to \E|_{C_n}\to 0.
	\]
	
	Since $h^0(S,\E(-nC))=h^1(S,\E(-nC))=0$, (the vanishing of $h^0$ follows from the defining sequence (\refeq{eqn:E})), we infer that $h^0(S,\E)=p+3$. This leads to a contradiction with (\refeq{eq:no-sections}).

	\medskip
	
	Finally, we note that $\E|_{C_n}$ contributes to $\mathrm{Cliff}(C_n)$. We compute $\gamma(\E|_{C_n})=\mu(\E|_{C_n})-h^0(C_n,\E|_{C_n})+2$. We have proved that $H^1(S,\E(-nC))=0$, and hence $h^0(C_n,\E|_{C_n})=h^0(S,\E)=p+3$, implying
	\[
	\gamma(\E|_{C_n})=\frac{n(2a+2p+1)}{2}-p-1<n(2a-2p-1)-2
	\]
	by the assumption $a\ge 3p+2$, which concludes the proof.
	\endproof

	\begin{remark}
		As mentioned in the preamble of \cite[Section 4]{Farkas-Ortega}, Mercat's conjecture holds for any curve of genus $g$ and gonality $k$ if $g>2(k-1)(k-2)$. In our case, note that, as $S$ contains no $(-2)$--curve, the results of \cite{Ciliberto-Pareschi} imply that $\mathrm{gon}(C_n)=n(2a-2p-1)$. Since $g(C_n)=2an^2+1$, the above inequality fails, even though both expressions are quadratic in $n$. The gonality of $C_n$ is small compared to the genus, and yet not sufficiently small to satisfy the conditions for Mercat's conjecture.
	\end{remark}

	The same strategy yields the following improvement of the main result of \cite{Sengupta}.
	
	\begin{thm}
		\label{thm:VeryGeneral}
		Let $S$ be a $K3$ surface with $\mathrm{Pic}(S)=\langle C\rangle$, where $C$ is a smooth curve of genus $g\ge 2$. Denote by $k=\left[\frac{g+3}{2}\right]$ and by $\mathcal{E}$ the Lazarsfeld-Mukai bundle associated to a $g^1_k$, $A$ on $C$. Let $n\ge 2$ and $C_n\in|nC|$ be a smooth curve. If either $n\ge 3$, or $n=2$ and $g\ge 9$, then $\mathcal{E}|_{C_n}$ is semistable with $\gamma(\mathcal{E}|_{C_n})<\mathrm{Cliff}(C_n)$ and thus it is a counterexample to Mercat's conjecture. Furthermore, if $n\ge 3$, then $\mathcal{E}|_{C_n}$ is stable.
	\end{thm}

	\proof
	We proceed along the lines of the proof of  Theorem \ref{thm:semistability}. Note that $g(C_n)=n^2(g-1)+1$ and $\mathrm{Cliff}(C_n)=2(n-1)(g-1)-2$.

	We first establish that $\mathcal{E}|_{C_n}$ is semistable, and it is stable if $n\ge 3$. Suppose $\mathcal{E}|_{C_n}$ is unstable and consider 
	\[
	0\to \mathcal{O}_{C_n}(B)\to \mathcal{E}|_{C_n}\to \mathcal{O}_{C_n}(C-B)\to 0
	\]
	a destabilizing sequence. 	If  $\mathcal{O}_{C_n}(C-B)\ne \mathcal{O}_{C_n}$, then it contributes to the Clifford index of $C_n$, and 
	\[
	\mathrm{Cliff}(\mathcal{O}_{C_n}(C-B)) < \mu(\mathcal{E}|_{C_n})-2=n(g-1)-2\le \mathrm{Cliff}(C_n)
	\]
	which leads to a contradiction. Note that, if $n\ge 3$, we have the stronger inequality $\mu(\mathcal{E}|_{C_n})-2 < \mathrm{Cliff}(C_n)$.
	
	Therefore, the destabilizing sequence is, in fact,
	\[
	0\to \mathcal{O}_{C_n}(C)\to \mathcal{E}|_{C_n}\to \mathcal{O}_{C_n}\to 0
	\]
	and it follows that $h^0(C_n,\mathcal{E}|_{C_n})\ge h^0(C_n, \mathcal{O}_{C_n}(C)) = h^0(S, \mathcal{O}_S(C))=g+1$.
	
	We prove that the restriction map  $H^0(S,\mathcal{E})\to H^0(C_n,\mathcal{E}|_{C_n})$ is an isomorphism. Since $h^0(S,\mathcal{E})=g-k+3$, this will be in contradiction with the inequality above. The vanishing of $H^0(S,\mathcal{E}(-nC))$ follows immediately from the sequence (\ref{eqn:E}), twisted with $\mathcal{O}_S(-nC)$. The surjectivity of the restriction map reduces to $H^1(S,\mathcal{E}(-nC))= 0$ which is equivalent to the vanishing of $H^1(S,\mathcal{E}^\vee(nC))$. Consider the defining sequence
	\[
	0\to \E^\vee\to H^0(A)\otimes\mathcal{O}_S\to A\to 0,
	\]
	twist it by $\mathcal{O}_S(nC)$ and take the long cohomology sequence. This reduces the problem to proving the surjectivity of the multiplication map 
	\[
	H^0(C,A)\otimes H^0(S,\mathcal{O}_S(nC))\to H^0(C,A(nC)).
	\]
	Since the restriction map $H^0(S,\mathcal{O}_S(nC))\to H^0(C,\mathcal{O}_C(nC))$ is surjective, it suffices to prove that the multiplication map
	\[
	H^0(C,A)\otimes H^0(C,\mathcal{O}_C(nC))\to H^0(C,A(nC))
	\]
	is surjective. To verify this, we apply once again Green's explicit $H^0$ Lemma, \cite[Corollary (4.e.4)]{Green}, as in the proof of Theorem \ref{thm:semistability}. We observe that the condition 
	$$
	\mathrm{deg}(A)+\mathrm{deg}(\mathcal{O}_C(C_n))\ge 4g+2
	$$ 
	is verified for any $n\ge 2$, by the hypothesis.
	
	To finish the proof, we compute $\gamma(\mathcal{E}|_{C_n})=\mu(\mathcal{E}|_{C_n})-h^0(\mathcal{E}|_{C_n})+2=(n-1)(g-1)+k-2<2(n-1)(g-1)-2$.
	\endproof

\end{document}